\magnification = \magstep1
\baselineskip =13pt
\overfullrule =0pt
\input amstex
\documentstyle{amsppt}
\vsize =500pt

\nologo
 
\topmatter

\title   Semiinfinite symmetric powers\endtitle
\rightheadtext{semiinfinite symmetric powers}

 \author{ M. Kapranov}\endauthor
\address  {Department of Mathematics,  University
of Toronto, 100 St. George St. Toronto On M5S 3G3 Canada
 }\endaddress

\email {  kapranov\@math.toronto.edu}\endemail

\endtopmatter

\document

\def\ind{{\lim\limits_{\longrightarrow}}}
\def\pro{{\lim\limits_{\longleftarrow}}}
\def\bd{{\big\downarrow}}
\def\0{{\overline{0}}}
\def\1{{\overline{1}}}

\noindent {\bf (0.1)} In this note  we develop a theory of measures, differential forms
and Fourier transforms on some infinite-dimensional
real vector spaces, by generalizing the following two
constructions:

\vskip .1cm

\noindent {\bf (a)} The construction of the semiinfinite wedge power of a space $V$
equipped with a class of commensurable subspaces [ADK] [KP] [S].
Recall that it is obtained as a certain double direct limit of
the exterior algebras of finite-dimensional subquotients of $V$. 

\vskip .1cm

\noindent {\bf (b)} The construction of the space of
 measures on a
nonarchimedean local field $K$ with maximal ideal
 {\bf m}
as a double inverse limit of the spaces of measures 
(which in this
 case are the same
as functions) on finite subquotients
 $\bold m^i/\bold m^j$ of $K$. 

\vskip .1cm

The importance of generalizing 
the construction (b) was emphasized by A.N. Parshin [P2-3].
Among other things,
 he pointed out that all the 4 types of possible
 combinations of inductive
and projective limits of the spaces of measures on 
$\bold m^i/\bold m^j$ have very transparent 
analytic meaning
(Schwartz-Bruhat functions vs. distributions,
 compact support
vs. arbitrary support) and called for a generalization of
this approach to higher
local fields. This note can be seen as step in that direction,
treating 2-dimensional local fields such as 
 $\Bbb R((t)), \Bbb C((t))$. More general
 2-dimensional local
fields and adeles (see [P1] for background) will be considered in a subsequent paper.

\vskip .3cm

\noindent {\bf (0.2)} Let us describe our constructions 
more
precisely. Our setting is exactly the same as 
the one needed in
the theory of the ``Japanese group $GL(\infty)$" and the
semiinfinite Grassmannians 
which encompasses (a). Though this is not universally known,
   the relevant class
of structured infinite-dimensional spaces was introduced in
the 1942 book of S. Lefschetz [L] under the name of
locally linearly compact spaces. We adopt this framework
which predates by 40 years the relatively recent interest in
``semiinfinite" structures. 

For a locally linearly compact $\Bbb R$-vector space $V$ we introduce
the space $M(V)$ of ``smooth measures" on $V$
as a double projective limit of the spaces of 
smooth measures
on the finite-dimensional subquotients $U_1/U_2$
 where $U_i$
run over  open, linearly compact subspaces of $V$.
So the procedure is quite similar to (b).
 Here, however, 
one encounters the  difficulty similar to that familiar in (a), namely
that smooth  measures can be restricted onto a subspace
only after tensoring with some 1-dimensional space of Haar measures,
and this leads to the fact that $M(V)$ will be only a projective
representation of the group $GL(V)$ of continuous automorphisms of $V$. 
This can be pinned down more precisely by saying that we have
a version $M_h(V)$ of $M(V)$ for each ``Haar theory" $h$ on $V$,
an object of a certain gerbe, see (2.3). 

The category $\Cal L$ of locally linearly compact spaces possesses
a perfect duality $V\mapsto V^\vee$ (the Lefschetz-Chevalley generalization
of Pontryagin duality)
and we construct a Fourier transform
from $M_h(V)$ to $M_{h^\vee}(V^\vee)$ where $h^\vee$ is a Haar theory
on $V^\vee$ naturally associated to $h$.

Further, we consider the de Rham complexes of forms on the subquotients
$U_1/U_2$ as before, and by taking their double projective
limit we construct the semiinfinite de Rham complex $\Omega^\bullet(V)$.
Its ``anomaly" (i.e., the central extension of $GL(V)$ arising from
the projective action) is much simpler than that of $M(V)$:
it reduces to a $\Bbb Z/2$-central extension coming from the
orientation, and in the case of $\Bbb C$-spaces vanishes altogether. 
The complex $\Omega^\bullet(V)$ is similar to 
$\Omega^{ch}$, the
chiral de Rham complex of Malikov, Schechtman and Vaintrob [MSV]. A  treatment of $\Omega^{ch}$ itself
from a similar standpoint  will be given
in the forthcoming paper [KV].

\vskip .3cm

\noindent {\bf (0.3)} The construction of (0.1)(a) can be
included in our framework as a particular case, if
we extend it to encompass super-vector spaces, see [M]. In particular,
for a purely odd  finite-dimensional super-vector space $W$
the space of measures  is the
exterior algebra of $W$, and the integration is understood
in the sense of Berezin. If we apply our approach to finite-dimensional
subqotients of a purely odd locally linearly compact supervector space $V$,
we get precisely the semiinfinite wedge space. Therefore
it is natural to call our space $M(V)$ (for any $V$)
the {\it semiinfinite symmetric power of $V$}. Unless we are in the purely
odd case, it is of a pronounced analytic flavor. If $V$ is equipped
with a positive definite quadratic form $q$, one can specify a more
algebraic subspace $AG(V,q)\subset M(V)$ invariant under the
orthogonal group $O(q)$ by considering measures which are
almost Gaussian, i.e., are product of a polynomial and a Gaussian measure.

Another natural candidate for a ``semiinfinite symmetric power"
of $V$ would be an irreducible module over the Heisenberg
algebra of $V\oplus V^\vee$. For example any open linearly compact subspace $U\subset V$
gives a ``vacuum module". But 
  unless $V$ is purely odd
(so that we have a Clifford algebra) these vacuum modules are
not isomorphic to each other, so there is no preferred one.
We show that our $M(V)$ (or, rather, the dual space $D(V)$
formed by distributions) contains naturally all such vacuum
modules. 

\vskip .3cm

\noindent {\bf (0.4)} It is perhaps worth emphasizing the difference of our
approach with the more traditional one
of probability theory and cylindric measures, cf. e.g., [K]. Namely, the latter amounts
to viewing an infinite-dimensional (say linear) space $V$ as a pro-object,
a certain completion of a projective limit of finite-dimensional spaces $W_i$.
 The approach of cylindric measures is (in modern language) an elaboration of
the naive idea of formally putting
$$\text{Meas}\, (\pro \,W_i) = \pro \, \text{Meas}\, (W_i),$$
where Meas stands for the space of measures.
 What we do is similar, except that we consider our spaces
not as pro-objects but rather as ind-pro-objects
which is the natural structure present on local fields. 

\vskip .3cm

\noindent {\bf (0.5)} This note was finished when
the author was staying at the Max-Planck-Institut
f\"ur Mathematik in summer 2001. The author's research
was also partially supported by grants from
NSF and NSERC.

\vfill\eject

\centerline {\bf \S 1. Locally linearly compact spaces}

\vskip 1cm

\noindent {\bf (1.1) Definitions.} We recall some classical concepts, due
to Lefschetz and Chevalley [L]. Let $k$ be a field, considered
with discrete topology. A topological $k$-vector space
is called {\it linearly topological} if it has a basis of neighborhoods
of 0  formed by linear subspaces. All linearly topological vector spaces
will be assumed Hausdorff. A linearly topological vector space $V$ is
called {\it linearly compact}, if any family $A_i\subset V$, $i\in I$
of closed affine subspaces such that $\bigcap_{i\in J} A_i\neq\emptyset$
for any finite $J\subset I$, has 
$\bigcap_{i\in I} A_i\neq\emptyset$. More generally, one says that $V$ is 
{\it locally linearly compact}, if it has a basis of neighborhoods of 0
formed by linearly compact open (automatically closed)
subspaces. We will drop ``linearly topological" when speaking about
(locally) linearly compact spaces. 

\vskip .1cm

\noindent {\bf (1.1.1) Examples.} (a) Any finite-dimensional $k$-vector
space with discrete topology is linearly compact. 
 The space $k[[t]]$ of formal Taylor series
with the $t$-adic topology is linearly compact. In fact, any product
of linearly compact spaces is linearly compact.

(b) The space $k((t))$ of formal Laurent series is locally linearly
compact. 

\vskip .2cm

\noindent We introduce the following categories:

\vskip .1cm

$\text{Vect}_0(k)$: finite-dimensional $k$-vector spaces.

\vskip .1cm

$\text{Vect}(k)$: all $k$-vector spaces (considered with discrete topology)

\vskip .1cm

$\widehat{\text{Vect}}(k)$: linearly compact $k$-spaces.

\vskip .1cm

$\Cal L(k)$: locally linearly compact spaces.

\vskip .1cm

\noindent Clearly, $\text{Vect}(k)$, $\widehat{\text{Vect}}(k)$ are full subcategories
of $\Cal L(k)$ with intersection $\text{Vect}_0(k)$. 

\vskip .1cm

\noindent {\bf (1.1.2) Remark.} From the modern
point of view, all the above categories can be obtained from
$\text{Vect}_0(k)$ by the purely algebraic constructions of passing to the
categories of ind- and pro-objects, see, e.g., [AM]. More precisely, in the
notations of {\it loc. cit.} $\text{Vect}(k)=\text{Ind}(\text{Vect}_0(k))$
and $\widehat{\text{Vect}}(k) = \text{Pro}(\text{Vect}_0(k))$. An embedding
in $\widehat{\text{Vect}}(k)$ is open iff it is induced (in the sense
of forming a Cartesian square)
from an embedding in $\text{Vect}_0(k)$. Finally, $\Cal L(k)$ is
identified with the full subcategory in 
$\text{Ind}(\text{Pro}(\text{Vect}_0(k)))$ formed by inductive systems 
over  $\text{Pro}(\text{Vect}_0(k))$ consisting of open (in the above sense)
embeddings. General definitions of locally compact ind-pro-objects (with
$\text{Vect}_0(k)$ replaced by more general categories) were proposed
by Beilinson [Be] and earlier by Kato [Ka]. We tried to keep this paper as
short and elementary as possible and thus avoided any systematic use
of ind/pro-objects. 

\vskip .3cm

\noindent {\bf (1.2) Dim and Det.} We start by  reformulating, in the form
needed for us, the basic constructions of [ADK]. Let $V$ be a locally
linearly compact $k$-vector space. We denote by $G(V)$ and call
the {\it semiinfinite Grassmannian of} $V$, the set of open
linearly compact subspaces in $V$. If $U_1, U_2\in G(V)$
and $U_1\subset U_2$, then $(U_2/U_1)$ is finite-dimensional.
Moreover, for arbitrary $U_1, U_2\in G(V)$
both $U_1\cap U_2$ and $U_1+U_2$ are in $G(V)$. 

By a {\it dimension theory} on $V$ we mean a map $d: G(V)\to \Bbb Z$ 
such that, whenever $U_1, U_2\in G(V)$
and $U_1\subset U_2$, we have $d(U_2)=d(U_1)+\dim(U_2/U_1)$. The set
of dimension theories will be denoted $\text{Dim}(V)$. The group
$\Bbb Z$ acts on $\text{Dim}(V)$ by adding constant functions. It is clear
from the above that this makes $\text{Dim}(V)$ into a $\Bbb Z$-torsor. 

For $W\in \text{Vect}_0(k)$ let $\det(W)$ be the top exterior power of $W$. 
Recall that for any short exact sequence
$$0\to W'\to W\to W''\to 0$$
in $\text{Vect}_0(k)$
we have a natural identification
$$m_{W'WW''}: \det(W')\otimes\det(W'')\to\det(W), \leqno (1.2.1)$$
and these identifications are associative for any filtration
$W_1\subset W_2\subset W$ of length 2. 

\proclaim{(1.2.2) Definition} Let $V$ be a locally
linearly compact $k$-vector space. A determinantal theory on $V$ is a
rule $\Delta$ which associates to each $U\in G(V)$ a 1-dimensional
$k$-vector space $\Delta(U)$, to each embedded pair $U_1\subset U_2$,
$U_i\in G(V)$, an isomorphism
$$\Delta_{U_1U_2}: \Delta(U_1)\otimes\det(U_2/U_1)\to \Delta(U_2)$$
so that for any nested triple $U_1\subset U_2\subset U_3$ the obvious
diagram
$$\matrix
&\Delta(U_1)\otimes\det(U_2/U_1)\otimes\det(U_3/U_3)&
 \longrightarrow&
\Delta(U_1)\otimes\det(U_3/U_1)&\\
&\bd&&\bd&\\
&\Delta(U_2)\otimes\det(U_3/U_2)& \longrightarrow &\Delta(U_3)&\\
\endmatrix$$
is commutative. 
\endproclaim

We denote by $\text{Det}(V)$ the category (groupoid) formed by all
determinantal theories on $V$ and their isomorphisms (in the obvious 
sense). If $\phi: \Delta\to \Delta'$ is an isomorphism of determinantal
theories and $\lambda\in k^*$, then $\lambda\phi$ is also an isomorphism.
One easily sees that:

\proclaim{(1.2.3) Proposition}
The above action of $k^*$ on the morphisms makes $\roman{Det}(V)$
into a $k^*$-gerbe [Bre] [Bry], i.e., each $\roman{Hom}_{\roman{Det}(V)}(\Delta,
\Delta')$ becomes a $k^*$-torsor and the composition is bilinear. 
\endproclaim

\vskip .1cm

\noindent {\bf (1.2.4) Example.} To compare our approach with that of [ADK],
note that any $U\in G(V)$ defines a unique dimension theory $d_U$
such that $d_U(U)=0$. The difference (in the $\Bbb Z$-torsor $\roman{Dim}(V)$)
of two such elements $d_U, d_{U'}$ is the integer denoted in [ADK]
by $[U':U]$ (the relative dimension of $U$ and $U'$). 

Similarly, a choice of $U$ gives rise to a canonical determinantal theory
$\Delta_U$ such that $\Delta_U(U)=k$ and the values of  $\Delta_U$
on other elements of $G(V)$ is recovered uniquely (that is, up to unique
isomorphism) from the axioms of Definition 1.2.2. The Hom-torsor
(in the $k^*$-gerbe $\text{Det}(V)$) between two such theories
$\Delta_U, \Delta_{U'}$ is the $k^*$-torsor corresponding to
the 1-dimensional vector space denoted in [ADK] by $[U'|U]$
(the relative determinant of $U'$ and $U$). 

\vskip .3cm

\noindent {\bf (1.3) The group $GL(\infty)$.} 
Let $V$ be a locally
linearly compact $k$-vector space. We write $GL(V)$ for $\text{Aut}_{
\Cal L(k)}(V)$, i.e., for the group of continuous automorphisms of $V$. 
When $V=k((t))$, this is the so-called ``Japanese group
$GL(\infty)$". 

The group $GL(V)$ acts on the $k^*$-gerbe $\text{Det}(V)$ and any object
$\Delta\in\text{Det}(V)$ gives, in a standard way, a central extension
$$1\to k^*\to \widetilde{GL}_\Delta(V)\to GL(V)\to 1.\leqno (1.3.1)$$
More precisely, see, e.g., [Bry], an element of 
$\widetilde{GL}_\Delta(V)$ is a pair $(g,\phi)$ where $g\in GL(V)$ and
$\phi: \Delta\to g(\Delta)$ is an isomorphism. The composition is
$(g_1, \phi_1) (g_2, \phi_2) = (g_1, g_2, g_2(\phi_1)\circ \phi_2)$.
Let us recall its standard properties.

\proclaim{(1.3.2) Proposition} (a) The extension (1.3.1) splits
iff $\Delta$ can be made into a $GL(V)$-equivariant object
of $\roman{Det}(V)$.

(b) If we consider the action of $\widetilde{GL}_\Delta(V)$
on $\roman{Det}(V)$ via the projection onto $GL(V)$, then $\Delta$
canonically has the structure of a $\widetilde{GL}_\Delta(V)$-equivariant
object.

(c) If $\Delta, \Delta'$ are two determinantal theories, then
$\roman{Hom}_{\roman{Det}(V)}(\Delta, \Delta')$ is identified with
$\roman{Hom}(\widetilde{GL}_\Delta(V), \widetilde{GL}_{\Delta'}(V)0$
in the category of central extensions. In parituclar, all the 
$\widetilde{GL}_\Delta(V)$ are isomorhic to each other. 

\endproclaim

\vskip .1cm

\noindent {\bf (1.3.3) Examples.} (a) A choice of a reference subspace $U\in
G(V)$ produces an object $\Delta_U\in\text{Det}(V)$ and hence a particular
central extension.

(b) If $V=k((t))$, $U=k[[t]]$, then $\widetilde{GL}_{\Delta_U}(V)$
  is the standard central
extension of the Japanese group [ADK]. 

(c) If $V$ is discrete, than $0\in G(V)$ so $\Delta_0$ is an equivariant
object of $\text{Det}(V)$ and the extension splits. Similarly, if $V$ is linearly
compact, then $V\in G(V)$ and the extension splits.

\vskip .3cm

\noindent {\bf (1.4) Further properties of Dim and Det.} The category
$\Cal L(k)$ possesses the following extra structures:

\vskip .1cm

\noindent (1.4.1) Duality, which is an antiequivalence $V\mapsto
V^\vee = \text{Hom}(V, k)$ (the space of continuous functionals).
The space $V^\vee$ is again locally linearly compact and
$G(V)$ and $G(V^\vee)$ are in order-reversing bijection $U\mapsto U^\perp$
(the orthogonal complement). 

\vskip .1cm

\noindent (1.4.2) The structure of an exact category in the sense of Quillen
[Q], i.e., a class of admissible short exact sequences
$$0\to V'\buildrel \alpha\over\longrightarrow V\buildrel\beta\over
\longrightarrow V''\to 0.\leqno (1.4.3)$$
More precisely, one calls a sequence (1.4.3)
admissible, if it is exact as a sequence of algebraic vector spaces,
if $V'$ is closed in $V$ and the topology on $V''$ is the quotient one.

These structures are compatible in the obvious sense: the dual of
an admissible short exact sequence is again admissible. As in [Q],
we will also speak about admissible filtrations etc. 

\vskip .1cm

On the other hand, for an abelian group $A$ the category of $A$-torsors
is a symmetric monoidal category with duality, see [Bre]. The monoidal
operation (tensor product of torsors over $A$) will be denoted $\otimes$
for $A=k^*$ and $\odot$ for $A=\Bbb Z$. The dual of a torsor $T$
will be denoted $T^\vee = \text{Hom}(T, A)$. For $t\in T$ let
$t^\vee\in T^\vee$ be the unique morphism taking $t$ to the
unit element of $A$. 

\proclaim {(1.4.4) Proposition} (a) For $V\in {\Cal L}(k)$ we have a canonical
identification of $\Bbb Z$-torsors $\roman{Dim}(V^\vee) = \roman{Dim}(V)^\vee$.

(b) For each admissible short exact sequence (1.4.3) we have a natural
identification of $\Bbb Z$-torsors 
$$\roman{Dim}(V')\odot \roman{Dim}(V'')\to \roman{Dim}(V)$$
and these identifications are associative in any admissible
filtration of length 2. 
\endproclaim

\noindent {\sl Proof:} (a) For a dimension theory $d$ on $V$
 we have a dimension
theory $d^\vee$ on $V^\vee$ given by $d^\vee(U) = -d(U^\perp)$, 
$U\in G(V^\vee)$.

(b) Given dimension theories $d'$ on $V'$ and $d''$ on $V''$,
we have a dimension theory $d$ on $V$ given by
$$d(U) = d'(\alpha^{-1}(U)) + d''(\beta(U)).$$
We leave the details to the reader. 

\vskip .2cm

For any $k^*$-gerbes $\frak G', \frak G''$ we denote $\frak G'\boxtimes
\frak G''$ the gerbe whose class of objects is $\text{Ob}(\frak G')\times
\text{Ob}(\frak G'')$ and
$$\text{Hom}_{\frak G'\boxtimes\frak G''}((x', x''), (y', y'')) = 
\text{Hom}_{\frak G'}(x', y') \otimes \text{Hom}_{\frak G''}(x'', y'').$$
We also denote by $\frak G^{op}$ the gerbe opposite to $\frak G$. We think of it as having for objects formal symbols $x^\vee, x\in\text{Ob}(\frak G)$
with $\text{Hom}_{\frak G^{op}}(x^\vee, y^\vee) = \text{Hom}_{\frak G}
(y,x)$. 

\proclaim{(1.4.5) Proposition} (a) For $V\in \Cal L(k)$ we have
a canonical equivalence of $k^*$-gerbes
$$\roman{Det}(V^\vee)\sim \roman{Det}(V)^{op}.$$
(b) for an admissible short exact sequence (1.4.3) we have a natural
equivalence of $k^*$-gerbes
$$\delta_{V'VV''}: \roman{Det}(V')\boxtimes\roman{Det}(V'')
\to\roman{Det}(V)$$
and natural transformation of ``associativity" for these equivalences
for any admissible filtration of $V_1\subset V_2\subset V$ of length 2:
$$\matrix
\roman{Det}(V_1)\boxtimes\roman{Det}(V_2/V_1)\boxtimes\roman{Det}(V/V_2)&
\longrightarrow &\roman{Det}(V_1)\otimes\roman{Det}(V/V_1)\\
\bd& \swarrow\varepsilon_{V_1\subset V_2\subset V}&\bd\\
\roman{Det}(V_2)\boxtimes\roman{Det}(V/V_2)&\longrightarrow &\roman{Det}(V)\\
\endmatrix$$
The transformations $\varepsilon_{V_1\subset V_2\subset V}$ fit into
a commutative cube for any admissible length 3 filtration.
\endproclaim

\noindent {\sl Proof:} (a) For a determinantal theory $\Delta$
on $V$ we have a determinantal theory $\Delta^\vee$ on $V^\vee$
given by $\Delta^\vee(U) = \Delta(U^\perp)^*$, $U\in G(V^*)$.

(b) (Sketch) Given determinantal theories $\Delta'$ on $V'$ and
$\Delta''$ on $V''$, we have a determinantal theory $\Delta = 
\delta_{V'VV''}(\Delta', \Delta'')$ on $V$ defined by
$$\Delta(U) = \Delta'(\alpha^{-1}(U))\otimes\Delta''(\beta(U)).$$
The (somewhat lengthy) checking of details as well as the
construction of the $\varepsilon_{V_1\subset V_2\subset V}$ and
verification of their properties, are left to the reader.

\vfill\eject

\centerline {\bf \S 2. Measures on locally linearly compact spaces}

\vskip 1cm

\noindent {\bf (2.0) Orientation issues.} From now on we take $k={\Bbb R}$
and write $\roman{Vect}_0$ for $\roman{Vect}_0(\Bbb R)$, as well as
$\Cal L$ for $\Cal L(\Bbb R)$ etc. Recall that 1-dimensional $\Bbb R$-vector
spaces are essentially the same as $\Bbb R^*$-torsors. If $L$ is such a space,
we denote by $|L|$ the 1-dimensional $\Bbb R$-vector space whose
corresponding torsor is induced from that of $L$ by the homomorphism
$$\Bbb R^*\to \Bbb R^*, \quad x\mapsto |x|.$$
Alternatively, $|L|$ can be identified with the space of functions
$s: L-\{0\}\to \Bbb R$ such that $s(\lambda x) = |\lambda|^{-1} s(x)$ for
any $\lambda\in \Bbb R^*$. 

Let  $W\in\text{Vect}_0$. Note that
the space $|\det(W)^*|$ is canonically identified
with the space of Haar measures on $W$. 
Further, let $\text{OR}(W)$ be the $\{\pm 1\}$-torsor of orientations
of $W$. Its two elements can be viewed as the two connected components
of the space $\Lambda^{\dim(W)}(W)-\{0\}$. Alternatively, they
can be viewed as the connected components of the space of all bases of $W$. 
Any $\{\pm 1\}$-torsor $O$ gives rise to a 1-dimensional $\Bbb R$-vector
space $O_{\Bbb R}$ via the canonical embedding $\{\pm 1\}\subset \Bbb R^*$. 
Explicitly, $O_{\Bbb R}$ can be viewed as consisting of 
odd (i.e., $\{\pm 1\}$-equivariant) functions
$s: O\to \Bbb R$. 

For any $W\in\text{Vect}_0$ we have a canonical identification
$$|\det(W)^*|\simeq \det(W)^*\otimes \text{OR}(W)_{\Bbb R}$$ 
which expresses the fact that a volume form in the presense of an
orientation gives a measure. 

\vskip .3cm

\noindent {\bf (2.1) The finite-dimensional case.}
 Let $W\in \roman{Vect}_0$
be a finite-dimensional $\Bbb R$-vector space. We introduce the following
function spaces: 

\vskip .2cm

$S(W)$: the Schwartz space of smooth rapidly decreasing functions
$W\to \Bbb R$, see [H].

\vskip .1cm

$D(W)$: the topological dual of $S(W)$, i.e., the space of
Schwartz distributions.

\vskip .1cm

$M(W) = S(W)\otimes|\det(W)^*|$: the space of smooth rapidly 
decreasing  measures on $W$.

\vskip .2cm

If $\beta: W\to W''$ is a surjection in $\text{Vect}_0$, we have the
direct image map (integration along the fibers):
$$\beta_*: M(W)\to M(W'').\leqno (2.1.1)$$
If $\alpha: W'\to W$ is an injection in $\text{ Vect}_0$,
 then the restriction of
functions induces a map
$$\alpha^*: M(W)\to M(W')\otimes|\det(\text{ Coker}(\alpha))^*|. \leqno (2.1.2)$$
The following is then straightforward.

\proclaim {(2.1.3) Proposition} (a) For two composable surjections $\beta_1,
\beta_2$ we have $(\beta_1\beta_2)_* = \beta_{1*}\beta_{2*}$.\hfill\break
(b) For two composable injections $\alpha_1, \alpha_2$
we have $(\alpha_1\alpha_2)_* = \alpha_{1*}\alpha_{2*}$, 
the equality understood with respect to the
canonical indetification $$\det(\roman{ Coker}(\alpha_1\alpha_2)) \simeq
\det(\roman{ Coker}(\alpha_1))\otimes \det(\roman{Coker}(\alpha_1))$$
(c) Let 
$$\matrix
&W&\buildrel {\alpha_2}\over\longrightarrow&W_1&\\
\beta_2&\big\downarrow &&\big\downarrow&\beta_1\\
&W_2&\buildrel {\alpha_1}\over\longrightarrow&W_{12}&\\
\endmatrix$$
be a  Cartesian square in $\roman{Vect}_0$ with $\beta_i$ being surjections
and $\alpha_i$ injections. (Such a square is authomatically cocartesian.) Then,
with respect to the identification $\roman{Coker}(\alpha_1)\simeq \roman{Coker}
(\alpha_2)$, we have the equality
$$\alpha_1^*\beta_{1*} = \beta_{2*}\alpha_1^*: M(W_1)\otimes |\det(\roman{Coker}
(\alpha_2))^*|\to M(W_2).$$
\endproclaim

In a similar way, an injection $\alpha: W'\to W$ defines   the direct image map
on distributions
$$\alpha_*: D(W')\to D(W),\leqno (2.1.4)$$
while a surjection $\beta: W\to W''$ gives rise to the inverse image map
$$\beta^*: D(W')\otimes|\det({\roman Ker}(\beta))^*|\to D(W). \leqno (2.1.5)$$
These maps, being dual to (2.1.1-2), satisfy the properties similar to
Propoisiton 2.1.3. 

\vskip .3cm

\noindent {\bf (2.2) Linearly compact case.} Let now $U$ be
a linearly compact space. We define the space of measures on $U$ as
$$M(U) = \pro_{U_1\subset U} M(U/U_1),\leqno (2.2.1)$$ 
where $U_1$ runs over open subspaces of $U$ (so that $U/U_1\in
\text{Vect}_0$) and the limit is taken with respect to the maps
(2.1.1) associated to the surjections $U/U_1\to U/U_2$, $U_2\subset U_1$. 
In a dual fashion, we define the space of distributions on $U$ to be
$$D(U) = \ind_{U_1\subset U} D(U/U_1)\otimes|\det(U/U_1)|,
\leqno (2.2.2)$$
where limit is now taken with respect to the maps obtained
by tensoring (2.1.5). 

\vskip .1cm

Let now $\alpha: U'\to U$ be an open embedding of linearly compact spaces.
For an open $U_1\subset U$ let $U_1'=\alpha^{-1}(U)$.
Then we have an embedding $\alpha_1: U'/U'_1\to U/U_1$ of finite-dimensional
spaces. If, further, $U_1\subset U_2\subset U$ are open and
$U'_2=\alpha^{-1}(U_2)$, then the square
$$\matrix
&U'/U'_1&\buildrel\alpha_1\over\longrightarrow&U/U_1&\\
\beta'&\bd&&\bd\beta\\
&U'/U'_2&\buildrel\alpha_2\over\longrightarrow &U/U_2&\\
\endmatrix$$
is Cartesian, so satisfies the assumptions of Proposition 2.1.3(c). We get,
therefore:

\proclaim {(2.2.3) Proposition} An open embedding
 $\alpha: U'\to U$ gives rise to the  inverse image map on measures
$$\alpha^*: M(U)\otimes|\det(\roman{Coker}(\alpha))^*|\to M(U')$$
and the direct image map on distributions
$$\alpha_*: D(U') \otimes|\det(\roman{Coker}(\alpha))^*|\to D(U),$$
and these maps are compatible with compositions.
\endproclaim

\vskip .2cm

\noindent{\bf (2.3) Locally linearly compact case.} Let now $V$ be a locally
linearly compact $\Bbb R$-vector space
and $G(V)$ the set of its open linearly compact subspaces.
 By a {\it Haar theory} on $V$ we
will understand a rule $h$ which to
each $U\in G(V)$ associates a 1-dimensional $\Bbb R$-vector space
$h(U)$ and to each pair $U_1\subset U_2$ an isomorphism
$h(U_1)\otimes |\det(U_2/U_1)^*|\to h(U_2)$ satisfying
a condition for nested triples similar to one given in
Definition 1.2.2. It is  clear that Haar theories form an $\Bbb R^*$-gerb,
in particular, any such theory $h$ gives rise to
a central extension $\widetilde{GL}_h(V)$ of $GL(V)$
by $\Bbb R^*$.
It is also clear that each determinantal theory
$\Delta$ on $V$ gives a Haar theory $|\Delta^*|$,
and $\widetilde{GL}_{|\Delta|^*}(V)$ is just the extension
induced from $\widetilde{GL}_\Delta(V)$ via the automorphism
$x\mapsto |x|^{-1}$ of $\Bbb R^*$.

\vskip .1cm

Fix a Haar theory
$h$ on $V$. We define the spaces of measures and distributions
on $V$ associated to $h$ to be
$$M_h(V) = \pro_{U\subset V} M(U)\otimes h(U)^* =
\pro_{U\subset V}\pro_{U_1\subset U} M(U/U_1)\otimes h(U)^*,
\leqno (2.3.1)$$
$$D_\Delta(V) = \ind_{U\subset V} D(U)\otimes h(U) = 
\ind_{U\subset V}\ind_{U_1\subset U} D(U/U_1)\otimes h(U).\leqno
(2.3.2)$$

\vskip .1cm

\noindent {\bf (2.3.3) Example.} If $V$ is discrete and $\Delta=\Delta_0$
is the determinantal theory described in (1.3.3)(c), then
$M_{|\Delta_0|^*}(V)$ is the space of Schwartz functions on $V$, 
i.e., the inverse limit of the Schwartz spaces on finite-dimensional
subspaces of $V$. Similarly,
if $V$ is linearly compact and $\Delta=\Delta_V$, then $M_{|\Delta_V|^*}(V) = 
M(V)$ is the inverse limit of the spaces of Schwartz measures on
finite-dimensional quotients of $V$. Thus in the general case,
$M_h(V)$ is a certain mixture of the spaces of functions and
measures.

\vskip .1cm

By construction, $M(V)$, $D(V)$ 
are representations of $\widetilde{GL}_h(V)$, formally dual to each other.
More precisely, we have a nondegenerate equivariant pairing
$$m, \phi\mapsto (m, \phi), \quad M(V)\otimes D(V)\to \Bbb R,
\leqno (2.3.4).$$
The evaluation of any particular $(m, \phi)$ reduces to 
the pairing of a function and a distribution on some finite-dimensional
subquotient of $V$. 

\vskip .3cm

\noindent {\bf (2.4) Heisenberg action on $M_h(V)$.}
Let $H(V)$ be the Heisenberg algebra generated by symbols $L_v, v\in V$,
$L_f, f\in V^*$, linearly depending on $v$ and $f$ and subject to the
relations
$$[L_f, L_v] =  f(v).$$
It can be viewed as the algebra of polynomial differential
operators on $V$. 

\proclaim{(2.4.1) Proposition} The space $M_h(V)$ has a natural
structure of a right $H(V)$-module, and $D_h(v)$ of a left $H(V)$-module.
\endproclaim

\noindent {\sl Proof:} Let $W$ be a finite-dimensional $\Bbb R$-vector
space and $H(W)$ be its Heisenberg algebra, defined as before.
Then, $H(W)$ is the algebra of polynomial linear differential
operators on $W$. As such, it acts on the left in $S(W)$ (functions)
and on the right
in $M(W)$ and $D(W)$ (measures and distributions). This would be 
the case for any smooth manifold. 
 
Let now $v\in V$ and
$U\subset V$ be any linearly compact
open subspace containing $v$. We define the operator $L_v$ on the
space $M(U) = \pro_{U'\subset U} M(U/U')$ to be induced by the morphism of
projective systems which on each $M(U/U')$ is given by the right action 
of the constant vector field
$v \, \text{mod}\, U'$. Then, we extend the action
to $M(U)\otimes h(U)^*$ by viewing $h(U)^*$ as a ``constant"
vector space of multiplicities. If $v\in U\subset U_1$, then
the resulting operators on $M(U)\otimes h(U)^*$ and
$M(U_1)\otimes h(U_1)^*$ are compatible
under the restriction map. Therefore we get an operator,
 still denoted $L_v$, on $M_h (V) = \pro_U M(U)\otimes
h(U)^*$.
The operators $L_f$ are defined in a similar way. Q.E.D.

\vskip .1cm

  For any linearly compact open subspace $U\i V$
we denote by $N_U$ the vacuum representation of $H(V)$
corresponding to $U$. By definition, it is generated by
one vector $|U\rangle$ (``vacuum") subject to the relations:
$$L_v |U\rangle =0, \quad v\in U, \quad\quad L_f |U\rangle = 0,\quad
f\in U^\perp.\leqno (2.4.2)$$
We can view it as the space of distributions on $V$ which
have support in $U$ and which are smooth along $U$. 
For different $U$ these modules are not isomorphic to each
other. In particular, the group $\widetilde{GL}(V)$ does not act, even projectively, in any
of the $N_U$ (even though its Lie algebra does, via
is embedding into a completion of $H(V)$ by vertex operators [KR]).

\proclaim {(2.4.3) Proposition}
 $N_U$ is naturally embeddeded into the space $D_{|\Delta_U|^*} (V)$
as an $H(V)$-submodule. Here $\Delta_U$ is the determinantal
theory associated to U, see (1.2.4). 
\endproclaim

\noindent {\sl Proof:} We exhibit an element $\delta_U$ (``delta
function along $U$") in $D_{|\Delta_U|^*}(V)$ satisfying the relations (2.4.2).
Then the embedding would be uniquely determined by sending 
 $|U\rangle$ to $\delta_U$. We represent $D_{|\Delta_U|^*}(V)$ as
$\ind_{U_1\supset U}\ind_{U_2\subset U_1} D(U_1/U_2)\otimes |\Delta_U(U_1)^*|$.
Consider the term of this double inductive system corresponding to
$U_1=U_2=U$. This term is canonically identified with $\Bbb R$
and we define $\delta_U$ as the image of $1\in\Bbb R$ in the
inductive limit. The relations (2.4.2) are verified straightforwardly. 
Q.E.D.

 \vskip .1cm

If $h$ is some other Haar theory, then we have
not a canonical embedding of $N_U$ into $D_h(V)$
but a canonical 1-dimensional space formed by such embeddings. 

\vskip .3cm

\noindent {\bf (2.5) Bilinear and quadratic forms.}
Let $V\in \Cal L$ be a locally linearly compact $\Bbb R$-space.
 We want to compare two possible approaches
to defining quadratic forms on $V$. First of all,
 by a bilinear form on $V$ we will mean
a morphism $b: V\to V^\vee$, where $V^\vee\in \Cal L$ is the 
topological dual of $V$.
 As usual, for such a $b$ we have
the transposed form $b^t: V\to V^*$ and $b$ is called symmetric
if $b^t=b$. We denote by $B(V)$ the space of all symmetric
bilinear forms on $V$. 
An element $b\in B(V)$ can be regarded as a continuous function
$b: V\times V\to \Bbb R$ and we denote by $q_b: V\to \Bbb R$
the associated quadratic form $q_b(x) = b(x,x)$. We call
$b$ nondegenerate if it is an isomorphism $V\to V^\vee$.
In this case we have a bilinear form $b^{-1}$ on $V^\vee$,
which is symmetric if $b$ is. We call $b$
 positive definite, if $q_b(x)>0$ for any $x\neq 0$. 
We denote by $B_{nd}(V)$ the set of all nondegenerate
symmetric bilinear forms on $V$ and by $B_+(V)\subset B_{nd}(V)$
the subset of forms which are both nondegenerate and positive definite.

Consider now the case of finite-dimensional vector spaces $W\in\text{Vect}_0$.
In this case the above concepts have their usual meaning. 
Let $\alpha: W'\to W$ be an injection in $\text{Vect}_0$
and $b\in B(W)$. Then we have the restriction $\alpha^*b\in B(W')$.
Note that for a nondegenerate $b$ the form $\alpha^*b$ may be degenerate,
but if $b$ is positive definite, then so is $\alpha^*b$
(and, in particular, $\alpha^*b$ is nondegenerate). 
Further,
  let $\beta: W\to W''$ be a surjection in $\text{Vect}_0$
and $b$ be a {\it positive definite} symmetric bilinear form
on $W$. Then we define the direct image $\beta_*(b)\in B(W'')$
as follows. Consider the injection $\beta^t: (W'')^*\to W^*$. 
dual to $\beta$ and define
$$\beta_*(b) = \bigl((\beta^t)^*(b^{-1})\bigr)^{-1}.\leqno (2.5.1)$$
Here we used the nondegeneracy of the positive
definite forms $b$ and $(\beta^t)^*(b^{-1})$. The following is then
elementary.

\proclaim{(2.5.2) Proposition} Let $q_b$ be the quadratic form
corresponding to $b$ and similarly for $\beta_*(b)$. Then
for $w''\in W''$ we have
$$q_{\beta_*b}(w'') = \max_{\beta(w)=w''} q_b(w'').$$
\endproclaim

\proclaim {(2.5.3) Corollary} (a) The operations  $\alpha^*$ and
$\beta_*$ on positive definitite symmetric bilinear forms 
(on finite-dimensional spaces) is compatible
with the composition of injections (resp. surjections).

(b)  Consider a Cartesian square in $\text{Vect}_0
(\Bbb R)$ as in (2.1.3)(c). Then, for any $q_1\in Q(W_1)$ we have $\alpha_1^*\beta_{1*}(q') = 
\beta_{2*}\alpha_1^*(q').$
\endproclaim

For $W\in \text{Vect}_0$ we denote by $Q(W)$ the
cone of positive definite quadratic forms on $W$. It is
canonically identified with $B_+(W)$. So we can think
of operations $\alpha_*$, $\beta_*$ as defined on elements of
$Q(W)$.

Let now $V\in\Cal L$ be a locally linearly compact $\Bbb R$-space. We  define
$$Q(V) = \pro_{U\subset V} \pro_{U'\subset U} Q(U/U'),\leqno(2.5.4)$$
where the limits are taken with respect to the operations
$\alpha^*$ and $\beta_*$ on quadratic forms. Elements of $Q(V)$
will be called (positive-definite) quadratic forms on $V$. 

\proclaim{(2.5.5) Proposition} The set $Q(V)$ is naturally
identified with $B_+(V)$.
\endproclaim

\noindent {\sl Proof:} Let $b$ be a symmetric bilinear form on $V$.
The topology in $V$ is that of inductive limit of its
linearly compact open subspaces $U\subset V$. Further, the
basis of linearly compact neighborhoods of 0 in $V^\vee$
is formed by the orthogonals $U^\perp$ where $U$ is a linearly
compact neighborhood of 0 in $V$. Accordingly, 
$$\text{Hom}_{\Cal L}(V, V^\vee) = 
\pro_{U_1\subset V}\ind_{U_2\subset V} \text{Hom} (U_1, U_2^\perp).$$
Note that in the above we can take $U_2\subset U_1$, and in this
case a morphism from $U_1$ to $U_2^\perp$ gives 
(after restriction to $U_1$) a bilinear
form $b_{U_1/U_2}$ on the finite-dimensional space $U_1/U_2$. Assuming that $b$
is nondegenerate and positive definite, we find that these
 forms on finite-dimensional subquotients are positive definite
and compatible with respect to the restrictions and projections.
In other words, the system of their associated quadratic forms 
is an element of the double inverse limit (2.5.4). 
The converse is similar. 

\vskip .3cm

\noindent {\bf (2.6) Gaussian measures.} Let $W\in\text{Vect}_0$.
Any $q\in Q(W)$. 
 $q\in Q(W)$ gives rise to a Haar measure
$\text{dVol}_q\in |\det (W^*)|$  (the measure of a $q$-orthocube is 1).  We get therefore the {\it  Gaussian measure}
$$ \gamma_q =  {1\over (2\pi)^{\dim(W)/2}} e^{-q(x)/2} \text{dVol}_q \in M(W).
\leqno (2.6.1)$$ 
 
Recall the standard properties of Gaussian integrals.

\proclaim{(2.6.2) Proposition} (a) If $\beta: W\to W''$ is a surjection
and $q \in Q(W)$, then
$\beta_*(\gamma_q)$ (the direct image of  measures as in (2.1.1))
is equal to $\gamma_{\beta_*(q)}$. 

(b) If $\alpha: W'\to W$ is an injection, $q\in Q(W)$ and $\roman
{dVol}_{\alpha, q}\in |\det(\roman{Coker}(\alpha))^*|$ is
the Haar measure induced by $q$ on the orthogonal complement of $\roman{Im}
(\alpha)$, then
$$\alpha^*(\gamma_q) = \gamma_{\alpha^*q}\otimes \roman{dVol}_{\alpha, q}.$$
\endproclaim

Let now $V\in \Cal L$ and $q\in Q(V)$. Note that $q$ trivializes the
space $|\det(U_2/U_1)^*|$ for any open, linearly comapct $U_2\subset U_2\subset
V$. This is because $U_2/U_1$ is identified with the $q$-orthogonal
complement of $U_1$ in $U_2$ and the latter comes equipped with the
Haar measure $\text{dVol}_q$. If $U_1\subset U_2\subset U_3\subset V$,
then these trivializations are compatible with the exact sequence
$$0\to U_2/U_1\to U_3/U_1\to U_3/U_2\to 0.$$
It follows that for any two open, linearly compact $U_1, U_2\subset V$
the Haar theories $|\Delta_{U_1}|^*, |\Delta_{U_2}|^*$ are canonically
identified. Further, Proposition 2.6.2 implies the following fact.

\proclaim{(2.6.3) Proposition} Fix $U\in G(V)$. Then any $q\in Q(V)$ gives rise to an element $\gamma_q\in M_{|\Delta_U|^*}(V)$ (the Gaussian measure)
invariant under the orthogonal group $O(V,q)$. 
\endproclaim

Let $W\in \text{Vect}_0$ and $q\in Q(W)$. A measure $\nu\in M(W)$ 
will be called {\it almost Gaussian}
(with respect to $q$) if it has the form $\nu = f\cdot\gamma_q$
where $f$ is a real polynomial function on $W$. We denote $AG(W,q)\subset M(W)$ the space of almost Gaussian measures. It is classical (``Wick's theorem")
that the class of almost Gaussian measures is preserved under the
operations of direct and inverse image, i.e., we have maps
$$\beta_*: AG(W,q)\to AG(W'', \beta_*q), \quad \beta: W\twoheadrightarrow W'',
\leqno (2.6.4)$$
$$\alpha^*: AG(W,q)\to AG(W', \alpha^*q)\otimes |\det(\text{Coker}(\alpha))^*|,
\quad \alpha: W'\hookrightarrow W.\leqno (2.6.5)$$

\proclaim {(2.6.6) Definition} Let $V\in \Cal L$ and $q\in Q(V)$. We
define $AG(V,q)$, the space of almost Gaussian measures on $V$ (with respect
to $q$) to be
$$AG(V,q) = \pro_{U\subset V} \pro_{U'\subset U} AG(U/U', q_{U/U'}),$$
where $q_{U/U'}$ is the quadratic form induced by $q$ on $U/U'$, see
(2.5.4) and the limits are taken with respect to the maps (2.6.4-5). 
\endproclaim

By construction, $AG(V,q)$ is a representation of the
group $O(q)$. It can be regarded as a kind of ``algebraic semiinfinite
symmetric power" of $V$.

\vskip .3cm

\noindent {\bf (2.7) The Fourier transform.}
Let $W\in\text{Vect}_0$. We denote by $S_{\Bbb C}(W)$, 
$M_{\Bbb C}(W)$ the complexifications of the spaces of Schwartz functions
and measures on $W$. Then, we have the Fourier transform which
is an isomorphism
$$\Cal F_W: M_{\Bbb C}(W)\to S_{\Bbb C}(W^*), \quad \Cal F(\mu)(y)
=\int_{x\in W} e^{i(x,y)}d\mu.\leqno (2.7.1)$$
Let  now $V\in\Cal L$ and $\Delta$ be a determinantal theory on $V$.
Let also $h^\vee$ be the dual determinantal theory on $V^\vee$,
defined as in (1.4.5).
 Then, the maps (2.7.1) for  finite-dimensional subquotients
of $V$ are assembled together into an isomorphism
$$\Cal F = \Cal F_V: M_h(V)\otimes \Bbb C\to M_{h^\vee}(V^\vee)
\otimes \Bbb C
\leqno (2.7.2)$$
which we will also call the Fourier transform. 
The following is an immediate consequence of the standard properties
of  the Fourier transform.

\proclaim{(2.7.3) Proposition} (a) The composition $\Cal F_{V^\vee}\circ
\Cal F_V$ is equal to the automorphism of $M_h(V)$
induced by the action of the element $(-1)\in GL(V)$. 

(b) If $q\in Q(V)$ and $q^{-1}\in Q(V^\vee)$ corresponds to the inverse
bilinear form, then $\Cal F_V(\gamma_q) = \gamma_{q^{-1}}$.

(c) In the situation of (b), the Fourier transform takes the
space $AG(V,q)\otimes\Bbb C$ into $AG(V^\vee, q^{-1})\otimes\Bbb C$. 

\endproclaim

\vfill\eject

\centerline{\bf \S 3. Forms on locally linearly compact spaces.}

\vskip 1cm

\noindent {\bf (3.1) The finite-dimensional case.} Let $W\in\text{Vect}_0$
be a finite-dimensional $\Bbb R$-vector space. Let
 $\Omega^\bullet(W)$ be the de Rham complex of differential forms
on $W$ whose components (with respect to some linear coordinate
system) are Schwartz functions. 
We denote
$$\widetilde{\Omega}^\bullet(W) = \Omega^\bullet(W)\otimes_{\Bbb R} 
\text{OR}(W)_{\Bbb R}.\leqno (3.1.1)$$
Then for each surjection $\beta: W\to W''$ in $\text{Vect}_0$
we have the direct image map (integration of forms along the
fibers in the presense of orientation)
$$\beta_*: \widetilde{\Omega}^\bullet(W)\to 
 \widetilde{\Omega}^\bullet(W'')[d], \quad d=\dim(\text{Ker}(\beta)).
\leqno (3.1.2)$$
Here $[d]$ means the shift of grading by $d$. Similarly, for an injection
$\alpha: W'\to W$ we have the restriction map
$$\alpha^*: \widetilde{\Omega}^\bullet(W)\to \widetilde{\Omega}^\bullet(W)
\otimes\text{OR}(\text{Coker}(\alpha))_{\Bbb R}.\leqno (3.1.3)$$

\vskip .1cm

\noindent {\bf (3.1.4) Remark.} If $W$ is a finite-dimensional $\Bbb C$-vector
space (considered as an $\Bbb R$-vector space) then $\text{OR}(W)$
is canonically trivialized. Accordingly, for surjections or
injections in $\text{Vect}_0(\Bbb C)$ we have functorialities
$\beta_*, \alpha^*$ sithout any twist (but with a shift of
grading for $\beta_*$). 

\vskip .1cm

\proclaim {(3.1.5) Proposition}
(a) Each $\beta_*$ and $\alpha^*$ is a morphism of complexes.

(b) The direct and inverse images of forms commute with compositions
of surjections (resp. injections).

(c) An analog of Proposition 2.1.3(c) holds. 

\endproclaim

\vskip .2cm

\noindent {\bf (3.2) The locally linearly compact case.} 
Let now $V\in\Cal L$ and $G(V)$, as in (1.2), denote the
set of open linearly compact subspaces in $V$. We
call an {\it orientation theory} on $V$ a rule $O$ which
to any $U\in G(V)$ associates a $\{\pm 1\}$-torsor
$O(U)$ and to any pair $U_1\subset U_2$ an isomorphism
$U(U_1)\otimes \text{OR}(U_2/U_1)\to O(U_2)$ so that for any
nested triple $U_1\subset U_2\subset U_3$ the diagram
analogous to (1.2.2) commutes. 

Every determinantal theory $\Delta$ on $V$ defines
an orientation theory $\text{OR}(\Delta)$ via the change
of structure groups given by $\text{sgn}: \Bbb R^*\to \{\pm 1\}$. 

\vskip .1cm

\noindent {\bf (3.2.1) Example.} Let $V\in\Cal L(\Bbb C)$ be a locally
linearly compact $\Bbb C$-space and $V_{\Bbb R}\in\Cal L$ be $V$
considered as an $\Bbb R$-space. Then any $U\in G(V_{\Bbb R})$ contains
a $\Bbb C$-subspace $U'\in G(V)$ so that $\dim_{\Bbb R}
(U/U')<\infty$. We set $C(U)=\text{OR}(U/U')$. A different'choice
 of $U'$ leads to a canonically isomorphic $\{\pm 1\}$-torsor
because a finite-dimensional $\Bbb C$-space has a canonical
orientation. It is clear that $C: U\mapsto C(U)$ is an
orientation theory on $V_{\Bbb R}$. 

\vskip .1cm

\proclaim {(3.2.2) Definition} Let $V\in\Cal L$ and $O$ be an
orientation theory on $V$. The (semiinfinite) de Rham complex
of $V$ associated to $O$ is
$$\Omega^\bullet_O(V) = \pro_{U\subset V}\pro_{U_1\subset U}
\widetilde{\Omega}^\bullet(U/U_1)\otimes\roman{OR}(U/U_1)_{\Bbb R}
[\dim(U/U_1)].$$
It is graded by the $\Bbb Z$-torsor $\roman{Dim}(V)$. 
\endproclaim

If $V\in \Cal L(\Bbb C)$ and $O=C$ is the canonical orientation
theory on $V_{\Bbb R}$ then we can write the above more simply
and denote it by
$$\Omega^\bullet(V) = \Omega^\bullet_C(V_{\Bbb R}) = 
 \pro_{U\subset V}\pro_{U_1\subset U} \Omega^\bullet((U/U_1)_{\Bbb R})
[2\dim_{\Bbb C}(U/U_1)],
\leqno (3.2.3)$$
where $U, U_1$ now run over open, linearly compact $\Bbb C$-subspaces
of $V$. This complex can be seen as an analog of the chiral
de Rham complex of Malikov, Schechtman and Vaintrob [MSV]. 

\vfill\eject

\centerline{\bf \S 4. Measures on locally linearly compact superspaces
and complexes}

\vskip 1cm

\noindent {\bf (4.1) Superspaces and Berezin integration.} Let $k$ be a field. 
By $\text{SVect}_0(k)$ we denote the category of finite-dimensional
$k$-supervector spaces [M]. Thus an object of 
$\text{SVect}_0(k)$ is a finite-dimensional $k$-vector space
$W$ together with a $\Bbb Z/2$-grading $W=W^\0 \oplus W^\1$; here
we denote elements of $\Bbb Z/2$ by $\0, \1$. By $\Pi$ we denote the functor
of shift of $\Bbb Z/2$-grading; $(\Pi W)^{\overline {i}} = W^{\overline
{i+1}}$.

Similarly, we denote by $S\Cal L(k)$ the category of locally
linearly compact $k$-supervector spaces $V=V^\0\oplus V^\1$. 
It can be identified with the category of
locally compact ind-pro-objects (in the sense of [Be] [K])
of $\text{SVect}_0(k)$.

\vskip .1cm

From now on we take $k=\Bbb R$ and write $\text{SVect}_0, S\Cal L$ etc. 
Let $W\in\text{SVect}_0$. According to the general principles
of super-analysis [M] we define $S(W)$, the Schwartz space of $W$ as
$$S(W)=S(W^\0)\otimes_{\Bbb R} \Lambda(W^{\1 *}),\leqno (4.1.1)$$
where $S(W^\0)$ is the usual Schwartz space and  $\Lambda(W^{\1 *})$
is the exterior algebra. similarly, the space of distributions on $W$
is defined as
$$D(W) = D(W^\0)\otimes\Lambda(W^\1).\leqno (4.1.2)$$
The 1-dimensional $\Bbb R$-space of
Haar measures on $W$ is defined to be
$$\mu(W) = |\det(W^\0)^*|\otimes\det(W^\1).\leqno (4.1.3)$$
This is justified by the existence of the integration map
$$\int_W: S(W)\otimes\mu(W)\to\Bbb R,\leqno (4.1.4)$$
which on $S(W^\0)\otimes|\det(W^\0)^*|$ is the usual integration
over $W^\0$ and on $\Lambda(W^{\1 *})\otimes\det(W^\1)$
is the ``Berezin integration" [M] which is, algebraically, just the projection
$$\Lambda(W^{\1 *})\otimes\det(W^\1) \to\Lambda^{\max}(W^{1 *}\otimes
\Lambda^{\max} (W^\1) \buildrel \sim\over\rightarrow \Bbb R.
\leqno (4.1.5)$$ 
We denote $M(W)=S(W)\otimes\mu(W)$ can call it, as in (2.1)
the space of (smooth, rapidly decreasing) measures on $W$. Generalizing (4.1.4),
for any surjection $\beta: W\to W''$ in $\text{SVect}_0$ we have the
direct image map
$$\beta_*: M(W)\to M(W''),\leqno (4.1.6)$$
and for any injection $\alpha: W'\to W$ the inverse image map
$$\alpha^*: M(W)\to M(W')\otimes\mu(\text{Coker}(\alpha)),
\leqno (4.1.7)$$
which satisfy the properties similar to those listed in Proposition 2.1.3. 

\vskip .3cm

\noindent {\bf (4.2) Measure on locally linearly compact superspaces.}
Let now $V=V^\0\oplus V^\1\in S\Cal L$. We denote by $G(V)$
the set of open linearly compact sub-superspaces $U\subset V$.
As every such $U$ is a direct sum $U=U^\0\oplus U^\1$, we have
that $G(V) = G(V^\0)\times G(V^\1)$. 
The concept of  a Haar theory from (2.3) generalizes trivially
to the case of a superspace such as $V$. Namely, a Haar theory on
$V$ is a rule $h$ which associates 
to any $U\in G(V)$ a 1-dimensional $\Bbb R$-vector space $h(U)$
and to any pair $U_1\subset U_2$ an isomorphism $h(U_1)\otimes\mu(U_2/U_1)
\to h(U_2)$ such that the properties of Definition 1.2.2 hold. 

Given a Haar theory $h$ on $V$, we define the spaces of measures and
distributions on $V$, similarly to (2.3), as
$$M_h(V) = \pro_{U\subset V} M(U)\otimes h(U)^* =
\pro_{U\subset V}\pro_{U_1\subset U} M(U/U_1)\otimes h(U)^*,
\leqno (4.2.1)$$
$$D_\Delta(V) = \ind_{U\subset V} D(U)\otimes h(U) = 
\ind_{U\subset V}\ind_{U_1\subset U} D(U/U_1)\otimes h(U),\leqno
(4.2.2)$$
where the limits are taken with respect to the maps (4.1.6-7) and  similar
maps on distributions. 

\vskip .1cm

\noindent {\bf (4.2.3) Example.} Let $V=V^\1$ be a purely odd superspace and
$\overline V = \Pi V$ be the same space but considered as an even one.
Then a Haar theory $h$ on $V$ is the same as a determinantal
theory $\Delta$ on $\overline V$. It follows that the space of
$h$-distributions on $V$ has the form
$$D_h(V) = \ind_{U\subset\overline{V}}\ind_{U_1\subset U}
\Lambda(U/U_1)\otimes\Delta(U) =: \Lambda^{{\infty\over 2}+\bullet}_\Delta
(\overline{V}).$$
In other words, this is the semiinfinite exterior power of
$\overline{V}$ corresponding to the determinantal theory $\Delta$, see 
[ADK] [KP] [S].
The space $M_h(V)$ is just the dual of $ \Lambda^{{\infty\over 2}+\bullet}_\Delta
(\overline{V})$ (being a double projective limit).
Note that $ \Lambda^{{\infty\over 2}+\bullet}_\Delta
(\overline{V})$ is graded by the $\Bbb Z$-torsor $\text{Dim}(\overline{V})$.

Therefore it is natural to think about $D_h(V)$ (even for a purely
even $V$) as a  ``semiinfinite symmetric power" of $V$. 

\vskip .3cm

\noindent {\bf (4.3) Measures on complexes.} Let $V^\bullet$ be a bounded
admissible complex over the exact category $\Cal L=\Cal L(\Bbb R)$.
By assembling together terms of the same parity, we associate to
$V^\bullet$ a superspace $V_{\sup}\in S\Cal L$ with
$$V^{\overline{i}}_{\sup} = \bigoplus_{j\equiv i \,\, (2)} V^j, \quad \overline
{I}\in\Bbb Z/2.\leqno (4.3.1)$$

A Haar theory on $V^\bullet$ is, by definition,
a Haar theory on $V_{\sup}$. For example, given determinantal theories
$\Delta_i$ on $V^i$ for each $i$, we have a Haar theory
$$h=\biggl(\bigotimes_{i\equiv 0 \, (2)} |\Delta_i|^* \biggr)
\otimes\biggl( \bigotimes_{i\equiv 1\, (2)} \Delta_i\biggr)\leqno (4.3.2)$$
on $V$. 

Let $h$ be a Haar theory on $V^\bullet$. We define the spaces of measures
and distributions on $V^\bullet$ to be those on $V_{\sup}$. 
These spaces are naturally graded by the $\Bbb Z/2$-torsor
$$\text{Dim}(V^\bullet)_{\sup} = \bigotimes_{i\equiv 1\, (2)} \text{Dim}(V^i)/
2\Bbb Z.\leqno (4.3.3)$$

\proclaim{(4.3.4) Proposition} The differential $d_V$ of $V^\bullet$
induces a natural differential $d$ in $M_h(V^\bullet)$, of
degree $\1$, satisfying $d^2=0$.
\endproclaim

\noindent {\sl Proof:} This can be seen as an instance of the
naturality of $M_h(V^\bullet)=M_h(V_{\sup})$ with respect
to the super-group $\widetilde{GL}_h(V_{\sup})$, if we regard
$d_V$ as an odd element of the Lie algebra of this group. More
precisely, consider the commutative $\Bbb R$-superalgebra
$\Lambda[\epsilon]$ generated by one odd  generator $\epsilon$
(with square 0). Then $V_{\sup}[\epsilon] = V_{\sup} \otimes_{\Bbb R}
\Lambda[\epsilon]$ is a locally linearly compact 
$\Lambda[\epsilon]$-supermodule.
we denote by $G_{\Lambda[\epsilon]}(V_{\sup}[\epsilon])$ the
set of linearly compact sub-$\Lambda[\epsilon]$-supermodules $U\subset
V_{\sup}[\epsilon]$ such that $V_{\sup}[\epsilon]/U$ is free
over $\Lambda[\epsilon]$. Then, we  repeat the definition
of Haar theories and the construction of $M_h$ in this setting
which is possible because of the functorial nature of
the Berezin integration.  Our statement is obtained
by considering the automorphism $1+\epsilon d_V$
of $V_{\sup}[\epsilon]$ and its action on $M_h$. 

\vskip .1cm

\noindent {\bf (4.3.5) Example: semiinfinite Koszul complex.}
 Let $V\in\Cal L$ and consider the complex
$CV=\{ V\buildrel \text{Id}\over\to V\}$ concentrated in degrees 0 and 1.
The gerbe of Haar theories of $CV$ can be written as $|\text{Det}(V)|\boxtimes
\text{Det}(V)^{op}$. This is equivalent to the gerbe of orientation theories
on $V$. Let such an orientation theory $O$ be chosen and $h$ be the
corresponding Haar theory on $CV$. The complex $M_h(CV)$ can be
called the {\it semiinfinite Koszul complex} of $V$. 
It is similar (but not identical) to
the semiinfinite de Rham complex $\Omega^\bullet_O(V)$.

\vfill\eject

\centerline {\bf References}

\vskip 1cm

\noindent [ADK] E. Arbarello, C. De Concini, V.G. Kac,
The infinite wedge representation and the reciprocity law for
algebraic curves, in:{\it  Proc. Symp. Pure Math.} {\bf 49} (1989),
part 1, p. 171-190. 

\vskip .1cm

\noindent [AM] M. Artin, B. Mazur, Etale Homotopy, 
Lecture Notes in Math. {\bf 100}, Springer-Verlag, 1969.

\vskip .1cm

\noindent [Be] A. A.  Beilinson, On the derived category of
perverse sheaves, in: ``K-theory, Arithmetic and Geometry"
(Y.I. Manin Ed.) (Lect. Notes in Math. {\bf 1289}), p. 27-41,
Springer-Verlag, 1987. 

\vskip .1cm

\noindent [Bre] L. Breen, On the Classification of 2-gerbes
and 2-stacks, {\it Ast\'erisque} {\bf 225}, Soc. Math. France, 1994. 

\vskip .1cm

\noindent [Bry] J.-L. Brylinski, Loop Spaces, Characteristic Classes
and Geometric Quantization (Progress in Math. {\bf 107}),
Birkh\"auser, Boston, 1993. 

\vskip .1cm

\noindent [H] L. H\"ormander, The Analysis of Linear Partial
Differential Operators I. Distribution Theory and Fourier Analysis
(Grund. math. Wiss. {\bf 256}) Springer-Verlag, 1990. 

\vskip .1cm

\noindent [Ka] K. Kato, The existence theorem for higher local
class field theory, math.AG/0012150
\vskip .1cm

\noindent [KP] V.G. Kac, D. Peterson, Spin and wedge representations
of infinite-dimensional Lie algebras and groups, {\it Proc. Nat. Acad. USA},
{\bf 78} (1981), 3308-3312. 

\vskip .1cm

\noindent [KR] V.G. Kac, A.K. Raina, Bombay Lectures on
highest-weight
representations of infinite-dimensional Lie algebras, World Scientific Publ.
Singapore, 1987. 

\vskip .1cm

\noindent [Ku] H.-H. Kuo, Gaussian Measures in Banach Spaces,
Lect. Notes in Math. {\bf 463}, Springer-Verlag, 1975.

\vskip .1cm

\noindent [KV] M. Kapranov, E. Vasserot, Vertex algebras and
the formal loop space, in preparation. 

\vskip .1cm

\noindent [L] S. Lefschetz, Algebraic Topology, (AMS Colloquium Publications
{\bf 27}), Amer. Math. Soc. New York, 1942.

\vskip .1cm

\noindent [MSV] F. Malikov, V. Schechtman, A. Vaintrob, Chiral de Rham
complex, {\it Comm. Math. Phys.} {\bf 204} (1999), 439-473.

\vskip .1cm

\noindent [M] Y.I. Manin, Gauge Fields and Complex Geometry,
Springer-Verlag, 1985. 

\vskip .1cm

\noindent [P1] A.N. Parshin, On the arithmetic of 2-dimensional schemes I.
Repartitions and residues, {\it Russian math. Izv.} {\bf 40} (1976),
736-773.

\noindent [P2] A.N. Parshin, Lectures at Universit\'e Paris-Nord,
June 1996. 

\vskip .1cm

\noindent [P3] A.N. Parshin, Higher-dimensional local
fields and L-functions, math.AG/0012151.  

\vskip .1cm

\noindent [Q] D. Quillen, Algebraic K-theory I, in: Lecture
Notes in Math. {\bf 341}, p. 85-147. 

\vskip .1cm

\noindent [S] M. Sato, The KP hierarchy and infinite-dimensional
Grassmann manifolds, 
in:{\it  Proc. Symp. Pure Math.} {\bf 49} (1989),
part 1, p. 51-66.

\bye